\documentclass[12pt]{article}

\usepackage{amssymb,amsmath}
\usepackage{theorem}
\usepackage{pifont}

\DeclareFontFamily{OT1}{nice}{}
\DeclareFontShape{OT1}{nice}{m}{n}{<5> <6> <7> <8> <9> <10>
<12><10.95><14.4><17.28><20.74><24.88>callig15}{}
\DeclareFontFamily{U}{nice}{} \DeclareFontShape{U}{nice}{m}{n}{<5>
<6> <7> <8> <9> <10>
<12><10.95><14.4><17.28><20.74><24.88>callig15}{}
\DeclareSymbolFont{calligra}{U}{nice}{m}{n}
\DeclareSymbolFontAlphabet{\nice}{calligra}
\DeclareFontFamily{OT1}{cmdh}{}
\DeclareFontShape{OT1}{cmdh}{m}{n}{<10>cmdunh10}{}
\input Umsa.fd
\input Umsb.fd
\input Ueuex.fd
\input Ueuf.fd
\input Ueur.fd
\input Ueus.fd
\def\nU{\mskip-4mu\nice U\mskip7mu{}}
\def\nV{\mskip-4mu\nice V\mskip10mu{}}

\def\nD{\nice D\mskip8mu{}}
\def\nE{\nice E\mskip8mu{}}
\def\nF{\nice F\mskip10mu{}}
\def\nG{\nice G\mskip8mu{}}
\def\nA{\nice A\mskip8mu{}}
\def\nB{\mskip-4mu\nice B\mskip7mu{}}
\def\nC{\nice C\mskip8mu{}}
\def\nH{\nice H\mskip10mu{}}

\newcommand {\art}[6]{{\sc #1:} {#2.} {\em #3} {\bf #4} {(#5),} {#6.}}

\newcommand {\book}[5]{{\sc #1:} {``#2."} {#3,} {#4} {(#5).}}
\newcommand {\samp}[8]{{\sc #1:} {#2,} {{\em in} ``#3,"} {(#4),}
                                  {#5,} {#6,} {(#7),} {#8.}}
\newcommand {\toappear}[3]{{\sc #1:} {#2.} {\em #3,} {\sf{to appear.}}}

\newcommand \reals {{[\omega]^{\omega}}}
\newcommand \fin {{[\omega]^{<\omega}}}
\newcommand \la {{\langle}}
\newcommand \ra {{\rangle}}
\newcommand \subs {{\subseteq}}
\newcommand \les {{\lesssim}}
\newcommand \surj {{\twoheadrightarrow}}
\newcommand \cf {{\sl cf.}\;}
\newcommand \eg {{\sl e.g.}}
\newcommand \ie {{\sl i.e.}}
\newcommand {\forces}[1] {{\hspace{0.5ex}}{\rule{0.1ex}{1.5ex}}
                    {\hspace{0.2ex}}{\rule{0.1ex}{1.5ex}}
                    {\rule[0.75ex]{1.5ex}{0.1ex}}
                    {\hspace{0.2ex}}{}_{#1}{\hspace{0.4ex}}}

\newcommand {\open}[1]{({#1})^{\omega}}
\newcommand {\nseg}[2]{({#1})^{#2*}}
\newcommand {\pf}[1]{\text{\tt{PF}}\big({#1}\big)}
\newcommand {\puf}[1]{\text{\tt{PUF}}\big({#1}\big)}
\newcommand \parti {(\omega)^{\le\omega}}
\newcommand \parto {(\omega)^\omega}
\newcommand \NN {(\N)}
\newcommand \ceq {\sqsubseteq}
\newcommand \kap {\sqcap}
\newcommand \seg {\preccurlyeq}

\newcommand \cP {{\mathcal{P}}}
\newcommand \cU {{\mathcal{U}}}
\newcommand \cI {{\mathcal{I}}}
\newcommand \cF {{\mathcal{F}}}
\newcommand \cA {{\mathcal{A}}}
\newcommand \cB {{\mathcal{B}}}

\newcommand \N {\mathbb{N}}
\newcommand \st {\flat}
\newcommand \MM {{\mathbb M}}
\newcommand \MMs {{\mathbb M}^{\mskip1mu\st}}
\newcommand {\MMUs}[1] {{\mathbb M}^{\mskip1mu\st}_{#1}}
\newcommand {\MMU}[1] {{\mathbb M}_{#1}}
\newcommand \UU {{\mathbb U}}
\newcommand \UUs {{\mathbb U}^{\mskip1mu\st}}

\newcommand \bfV {\mathbf{V}}
\newcommand \bfU {\mathbf{\check{U}}}
\newcommand \bfN {\mathbf{N}}

\newcommand \fp {{\mathfrak{p}}}
\newcommand \fri {{\mathfrak{i}}}
\newcommand \fb {{\mathfrak{b}}}
\newcommand \fs {{\mathfrak{s}}}
\newcommand \fc {{\mathfrak{c}}}

\newcommand \fu {{\mathfrak{u}}}
\newcommand \fU {{\mathfrak{U}}}
\newcommand \fh {{\mathfrak{h}}}
\newcommand \fH {{\mathfrak{H}}}
\newcommand \fr {{\mathfrak{r}}}
\newcommand \fR {{\mathfrak{R}}}
\newcommand \fhom {{\mathfrak{hom}}}
\newcommand \fHom {{\mathfrak{Hom}}}
\newcommand \fpar {{\mathfrak{par}}}
\newcommand \fPar {{\mathfrak{Par}}}
\newcommand \twotoc {{\mathfrak{2}}^\fc}

\newcommand \bN {\beta \N}
\newcommand \bNm {\beta \N\setminus \N}

\newcommand \add {\operatorname{add}}
\newcommand \cov {\operatorname{cov}}
\newcommand \non {\operatorname{non}}
\newcommand \fil {\operatorname{fil}}
\newcommand \mdom {\operatorname{dom}}
\newcommand \mmin {\operatorname{min}}
\newcommand \mmax {\operatorname{max}}
\newcommand \mMin {\operatorname{Min}}

\newcommand \ZFC {{\text{\sf ZFC}}}
\newcommand \MA {{\text{\sf MA}}}
\newcommand \CH {{\text{\sf CH}}}

\def\phi{\varphi}
\def\epsilon{\varepsilon}

\theorembodyfont{\rmfamily}
\newtheorem {nummersec}{ }[section]
\newtheorem {thmm}[nummersec]{{\sc{Theorem}}\rm\small}
\newtheorem {nummer}{ }[subsection]
\newtheorem {thm}[nummer]{{\sc{Theorem}}\rm\small}
\newtheorem {prop}[nummer]{{\sc{Proposition}}\rm\small}
\newtheorem {lm}[nummer]{{\sc{Lemma}}\rm\small}
\newtheorem {fct}[nummer]{{\sc{Fact}}\rm\small}
\newtheorem {cor}[nummer]{{\sc{Corollary}}\rm\small}

\newcommand \proof {\noindent{\sc Proof:}\hspace*{3mm}}

\def\eop{{\unskip\nobreak\hfil\penalty50\hskip8mm\hbox{}
  \nobreak\hfil
  {$\boldsymbol{\dashv}$}\parfillskip=0mm \par\smallskip}}

\begin{document}

\begin{center}
      {\Large{\bf{RAMSEYAN ULTRAFILTERS}}}
\end{center}
\smallskip

\begin{center}
{\sc{Lorenz Halbeisen\footnote{I would like to thank the {\it Swiss
National Science Foundation\/} for its support during the period in
which the research for this paper has been done.}}}\\[1.7ex]
{\small{\sl Department of Pure Mathematics\\ Queen's University
Belfast\\ Belfast BT7 1NN, Northern Ireland}}\\[1ex] {\small{\sl
Email: halbeis@qub.ac.uk}} \vspace{.5cm}
\end{center}

\begin{abstract}\noindent
We investigate families of partitions of $\omega$ which are related
to special coideals, so-called happy families, and give a dual form
of Ramsey ultrafilters in terms of partitions. The combinatorial
properties of these partition-ultrafilters, which we call Ramseyan
ultrafilters, are similar to those of Ramsey ultrafilters. For
example it will be shown that dual Mathias forcing restricted to a
Ramseyan ultrafilter has the same features as Mathias forcing
restricted to a Ramsey ultrafilter. Further we introduce an ordering
on the set of partition-filters and consider the dual form of some
cardinal characteristics of the continuum.
\end{abstract}
\renewcommand{\thefootnote}{}
\footnotetext{{}\hfill\\[-1.5ex] {\it 2000 Mathematics Subject
Classification:} {\bf 05D05 05D10} 03E05 03E40 03E17 03E35}

\setcounter{section}{-1}
\section{Introduction}

The Stone-\v{C}ech compactification $\bN$ of the natural numbers, or
equivalently, the ultrafilters over $\omega$, is a well-studied space
({\cf}{\eg}\;\cite{vMill} and \cite{ComfortNegrepontis}) which has a
lot of interesting topological and combinatorial features
({\cf}\cite{HindmanStrauss} and \cite{ToTopics}). In the late 1960's,
a partial ordering on the non-principal ultrafilters $\bNm$, the
so-called Rudin-Keisler ordering, was established and ``small''
points with respect to this ordering were investigated rigorously
({\cf}\cite{Booth}, \cite{BlassTrans}, \cite{BlassJSL} and
\cite{Laflamme}). The minimal points have a nice combinatorial
characterization which is related to Ramsey's Theorem
({\cf}\cite[Theorem\,A]{Ramsey}) and so, the ultrafilters which are
minimal with respect to the Rudin-Keisler ordering are also called
Ramsey ultrafilters (for further characterizations of Ramsey
ultrafilters see \cite[Chapter\,4.5]{BartoszynskiJudah}). Families,
not necessarily filters, having similar combinatorial properties as
Ramsey ultrafilters, are the so-called happy families
({\cf}\cite{Mathias}), which are very important in the investigation
of Mathias forcing ({\cf}\cite{Mathias}).

From the category theoretical point of view, subsets of $\omega$ and
partitions of $\omega$ are dual to each other (see
{\eg}\;\cite[Introduction]{luxor}), and therefore, it is natural to
look for the dualization of statements about subsets of $\omega$ in
terms of partitions of $\omega$. In this dualization process, a lot
of work is already done. Confer: \cite{luxor} for a dualization of
$\bN$; \cite{CarlsonSimpson}, \cite{Lorisym} and \cite{boise} for the
dualization of the Ramsey property and of Mathias forcing;
\cite{CarlsonSimpson} for a dualization of Ramsey's Theorem;
\cite{Cichon.etal} and \cite{Lorisha} for the dualization of some
cardinal characteristics of the continuum.

To investigate partition-filters, a useful tool is missing: the
dualization of Ramsey ultrafilters. The aim of this paper is to fill
this gap.

\section{Partition-filters}

\subsection{Notations and definitions}

Most of our set-theoretic notation is standard and can be found in
textbooks like \cite{Jechbook}, \cite{Kunen} or
\cite{BartoszynskiJudah}. So, we consider a natural number $n$ as an
ordinal, in particular $n=\{k:k<n\}$ and $0=\emptyset$, and
consequently, the set of natural numbers is denoted by $\omega$. For
a set $S$, $\cP(S)$ denotes the power-set of $S$. The notation
concerning partitions is not yet standardized. However, we will use
the notation introduced in \cite{Lorisym}.

A {\bf partition} $X$ of a set $S$ consisting of pairwise disjoint,
non-empty sets, such that $\bigcup X = S$. The elements of a
partition are called {\bf blocks}. Mostly, we will consider
partitions of $\omega$, so, if not specified otherwise, the word
``partition'' refers to a {\it partition of\/} $\omega$.

Most of the partitions in consideration are infinite, or in other
words, contain infinitely many blocks. However, at some places we
also have to consider finite partitions, this means, partitions
containing only finitely many blocks. The unique partition containing
just one block is denoted by $\{\omega\}$. The set of all partitions
is denoted by $\parti$ and the set of all partitions containing
infinitely many blocks is denoted by $\parto$.

Let $X$ and $Y$ be two partitions of a set $S$. We say $X$ is {\bf
coarser} than $Y$, or that $Y$ is {\bf finer} than $X$ (and write
$X\ceq Y$), if each block of $X$ is the union of blocks of $Y$. Let
$X\kap Y$ denote the finest partition of $S$ which is coarser than
$X$ and $Y$.

Further, for $n\in\omega$ and a partition $X\in\parti$, let
$X\kap\{n\}$ be the partition we get, if we glue all blocks of $X$
together which contain a member of $n$. If $X$ and $Y$ are two
partitions, then we write $X\ceq^* Y$ if there is an $n\in\omega$
such that $(X\kap \{n\})\ceq Y$.

A set $\nF\subs\parti$ is a {\bf partition-filter}, if the following
holds:
\begin{enumerate}
\item[{(a)}] $\{\omega\}\notin\nF$.
\item[{(b)}] For any $X,Y\in\nF$ we have $X\kap Y\in\nF$.
\item[{(c)}] If $X\in\nF$ and $X\ceq Y\in\parti$, then $Y\in\nF$.
\end{enumerate}

A partition-filter $\nF\subs\parti$ is called {\bf principal}, if
there is a partition $X\in\parti$ such that $\nF=\{Y:X\ceq Y\}$.

A set $\nU\subs\parti$ is a {\bf partition-ultrafilter}, if $\nU$ is
a partition-filter which is not properly contained in any
partition-filter.

Notice that a partition-ultrafilter $\nU$ which does not contain a
finite partition is always non-principal, and vice versa, a principal
partition-ultrafilter always contains a finite partition, in fact it
contains a $2$-block partition (see \cite[Fact\,3.1]{luxor}). Thus,
if $\nU$ is a non-principal partition-ultrafilter, $X\in\nU$ and
$X\ceq^* Y$, then $Y\in\nU$.

In the sequel we are mostly interested in partition-filters which
don't contain a finite partition, or in other words, in
partition-filters $\nF\subs\parto$.

For the sake of convenience, we defined the notion of
partition-filter only for partition-filters over $\omega$, but it is
obvious how to generalize this notion for partition-filters over
arbitrary sets $S$ (see also \cite{luxor}).

\subsection{An ordering on the set of
partition-filters}\label{sec:ordering}

Let $\pf{\parti}$ denote the set of all partition-filters. We define
a partial ordering on $\pf{\parti}$ which has some similarities with
the Rudin-Keisler ordering on $\bNm$.

To keep the notation short, for $\nH\subs\cP\big(\cP(\omega)\big)$
and a function $f:\omega\to\omega$ we define
$$f^{-1}(\nH):=\{f^{-1}(X):X\in\nH\}\,,$$ where for $X\in\nH$ we
define $$f^{-1}(X):=\{f^{-1}(b):b\in X\}\,,$$ where for
$b\subs\omega$, $f^{-1}(b):=\{n:f(n)\in b\}$.

Let $f:\omega\surj\omega$ be any surjection from $\omega$ onto
$\omega$ and let $X\in\parti$ be any partition. Then $f(X)$ denotes
the finest partition such that whenever $n$ and $m$ lie in the same
block of $X$, then $f(n)$ and $f(m)$ lie in the same block of $f(X)$.

For any partition-filter $\nF\in\pf{\parti}$ define
$$f(\nF):=\big{\{}Y\in\parti:\exists X\in\nF\big(f(X)\ceq
Y\big)\big{\}}\,.$$

We define the ordering ``$\les$'' on $\pf{\parti}$ as follows:
$$\nF\les\nG\text{\ \ if and only if\ \ }\nF =f(\nG)\text{ for some
surjection $f:\omega\surj\omega$}\,.$$

Since the identity map is a surjection and the composition of two
surjections is again a surjection, the partial ordering ``$\les$'' is
reflexive and transitive.

\begin{fct}\label{fct:supfilter}
Let $\nF,\nG\in\pf{\parti}$ and assume $f(\nG)=\nF$ for some
surjection $f:\omega\surj\omega$. Then $\nG\subs f^{-1}(\nF)$ and
$f^{-1}(\nF)\in\pf{\parti}$.
\end{fct}

\proof Let $\nH=f^{-1}(\nF)$, where $f:\omega\surj\omega$ is such
that $f(\nG)=\nF$. Since $\nF$ is a partition-filter and $f$ is a
function, for any $X_1,X_2\in\nF$ we have $X_1\kap X_2\in\nF$ and
$f^{-1}(X_1\kap X_2)=f^{-1}(X_1)\kap f^{-1}(X_2)$, and therefore,
$\nH$ is a partition-filter. Further, for any $Y\in\nG$ we get
$f(Y)\in\nF$ and $f^{-1}(f(Y))\ceq Y$, which implies $\nG\subs \nH$.
\eop

The ordering ``$\les$'' induces in a natural way an equivalence
relation ``$\simeq$'' on the set of partition-filters $\pf{\parti}$:

$$\nF\simeq\nG\text{\ \ if and only if\ \ }\nF\les\nG\text{ and
}\nG\les\nF\,.$$

So, the ordering ``$\les$'' induces a partial ordering of the set of
equivalence classes of partition-filters. Concerning
partition-ultrafilters, we get the following.

\begin{fct}\label{fct:isomorph}
Let $\nU,\nV\in\puf{\parti}$ and assume that $\nU$ is principal or
contains a partition, all of whose blocks are infinite. If
$\nU\simeq\nV$, then there is a permutation $h$ of $\omega$ such that
$h(\nU)=\nV$.
\end{fct}

\proof Because $\nU\les\nV$ and $\nV\les\nU$, there are surjections
$f$ and $g$ from $\omega$ onto $\omega$ such that $\nV =f(\nU)$ and
$\nU =g(\nV)$, and because $\nU$ and $\nV$ are both
partition-ultrafilters, by Fact~\ref{fct:supfilter} we get $\nU
=f^{-1}(\nV)$ and $\nV =g^{-1}(\nU)$.\\ First assume that $\nU$ is
principal and therefore contains a $2$-block partition
$X=\{b_0,b_1\}$. Because $g^{-1}(X)\in\nV$, the partition-ultrafilter
$\nV$ is also principal and we get $\nV=\{Y\in\parti:g^{-1}(X)\ceq
Y\}$, where $g^{-1}(X)=\{g^{-1}(b_0),g^{-1}(b_1)\}=:\{c_0,c_1\}$.
Now, because $\nU =f^{-1}(\nV)$, we must have $f^{-1}\big(
g^{-1}(X)\big)=X$, which implies $f^{-1}\big(g^{-1}(b_i)\big)
\in\{b_0,b_1\}$ (for $i\in\{0,1\}$). If one of the blocks of $X$ is
finite, say $b_0$, then $f|_{b_0}$ as well as $g|_{f(b_0)}$ must be
one-to-one, and therefore, $b_0$ has the same cardinality as $c_0$.
Hence, no matter if one of the blocks of $X$ is finite or not, we can
define a permutation $h$ of $\omega$ such that $h(b_0)=c_0$ and
$h(b_1)=c_1$, which implies $h(\nU)=\nV$.\\ Now assume that $\nU$
contains a partition $X=\{b_i:i\in\omega\}$, all of whose blocks
$b_i$ are infinite. Because $g$ is a surjection, $g^{-1}(X)$, which
is a member of $\nV$, is a partition, all of whose blocks are
infinite. Let $h$ be a permutation of $\omega$ such that
$h(b_i)=g^{-1}(b_i)$. Take any $Y\in\nV$ with $Y\ceq g^{-1}(X)$. By
the definition of $h$ we have $h^{-1}(Y)=g(Y)$ and since $\nU
=g(\nV)$ there is a $Z\in\nU$ such that $g(Y)=Z$, which implies
$h(Z)=Y$, hence, $h(\nU)=\nV$. \eop

The following proposition shows that ``$\les$'' is directed upward
(for a similar result concerning the Rudin-Keisler ordering see
\cite[p.\,147]{BlassTrans}).

\begin{fct}\label{fct:iso}
For any partition-filters $\nD,\nE\in\pf{\parti}$, there is a
partition-filter $\nF\in\pf{\parti}$, such that $\nD\les\nF$ and
$\nE\les\nF$.
\end{fct}

\proof Let $\rho_1$ and $\rho_2$ be two functions from $\omega$ into
$\omega$ defined by $\rho_1(n):=2n$ and $\rho_2(n):=2n+1$. For a
partition $X$ and $i\in\{0,1\}$, let $\rho_i(X):=\{\rho_i(b):b\in
X\}$, where $\rho_i(b):=\{\rho_i(n):n\in b\}$. Now, take any two
partition-filters $\nD,\nE\in\pf{\parti}$ and define $\nF$ by $$\nF
:=\big{\{}\rho_1(X)\cup \rho_2(Y): X\in\nD\wedge Y\in\nE\big{\}}\,.$$
Clearly, this defines a partition-filter. Define two surjections $f$
and $g$ from $\omega$ onto $\omega$ as follows:
\begin{equation*} f(n)=\begin{cases} \frac{n}{2}& \text{if
$n$ is even,}\\ 0& \text{otherwise.}\end{cases}\end{equation*}
\begin{equation*} g(n)=\begin{cases} \frac{n-1}{2}& \text{if
$n$ is odd,}\\ 0& \text{otherwise.}\end{cases}\end{equation*} It is
easy to verify that $f(\nF)=\nD$ and $g(\nF)=\nE$, which implies
$\nD\les\nF$ and $\nE\les\nF$. \eop

\section{Ramseyan ultrafilters}

\subsection{Coloring segments}\label{sec:coloring}

If $X$ is a partition of a set $S$, then we say that $S$ is the {\bf
domain} of $X$, written $\mdom(X)=S$. The set of all partitions of
natural numbers $n\in\omega$, called {\bf segments}, is denoted by
$\NN$. Thus, $s\in\NN$ implies $\mdom(s)\in\omega$. In particular,
$\emptyset$ is the unique partition of $0$ and
$\{\{\emptyset\}\}=\{1\}$ is the unique partition of $1$. For
$s\in\NN$, $|s|$ denotes the cardinality of $s$, which simply means
the number of blocks of $s$, and $\bigcup s:=\{\mdom (s)\}$.

For a set $b\subs\omega$, let $\mmin(b)$ be the least element of $b$
and for a set $P\subs\cP(\omega)$, let $\mMin(P):=\{\mmin(b): b\in P
\}$. Further, for a finite set $b\subs\omega$, let $\mmax(b)$ be the
greatest element of $b$. For $X\in\parti$, $s\in\NN$ and
$n\in\omega$, let $X(n)$ and $s(n)$ be the $n$th block of $X$ and
$s$, respectively, where we start counting with $0$ and assume that
the blocks are ordered by their least element.

Let $s,t\in\NN$ and $X\in\parti$: We write $s\ceq X$, if each block
$b\in s$ is the union of some sets $b_i\cap\mdom (s)$, where each
$b_i$ is a block of $X$; we write $s\seg t$ and $s\seg X$,
respectively, if for each $b\in s$ there is a $c_b\in t$ and a
$d_b\in X$, respectively, such that $b=c_b\cap\mdom (s)=d_b\cap\mdom
(s)$ (notice that $s\seg t$ implies $\mdom (s)\subs\mdom (t)$); and
for $s\ceq X$, $s\kap X$ denotes the finest partition $Y\in\parti$,
such that $s\seg Y\ceq X$.

For $s\in\NN$, let $s^*$ denote the partition $s\cup\big{\{}\{\mdom
(s)\}\big{\}}$. In particular, $\emptyset^* =\{1\}$. Notice that
$|s^*|=|s|+1$.

For $s\in\NN$ and $X\in\parto$ with $s\ceq X$, let $$\open{s,X}:=\{
Y\in\parto : s\seg Y\ceq X\}\,.$$ A set $\open{s,X}$, where $s$ and
$X$ are as above, is called a {\bf dual Ellentuck neighborhood}
({\cf}\cite[p.\,275]{CarlsonSimpson}). In particular,
$\open{\emptyset,X}=\open{\{1\},X}=: \open{X}$.

For $n\in\omega$, $\nseg{\omega}{n}$ denotes the set of all $u\in\NN$
such that $|u|=n$. Further, for $n\in\omega$ and $X\in\parto$ let
$$\nseg{X}{n}:=\big{\{}u\in\NN
:|u|=n\wedge u^*\ceq X\big{\}}\,;$$ and if $s\in\NN$ is such that $|s|\le n$
and $s\ceq X$, let $$\nseg{s,X}{n}:=\big{\{}u\in\NN
:|u|=n\wedge s\seg u\wedge u^*\ceq X\big{\}}\,.$$

{}From the so-called Dual Ramsey Theorem of Carlson and Simpson,
which is Theorem\,1.2 of \cite{CarlsonSimpson}, we get the following.

\begin{prop}\label{prop:thm}
For any coloring of $\nseg{\omega}{(n+1)}$ with $r+1$ colors, where
$r,n\in\omega$, and for any $Z\in\parto$, there is an infinite
partition $X\in\open{Z}$ such that $\nseg{X}{(n+1)}$ is
monochromatic.
\end{prop}

This combinatorial result is the dualization of Ramsey's Theorem,
\cite[Theorem\,A]{Ramsey}, in terms of partitions.

We say that a surjection $f:\omega\surj\omega$ {\bf respects} the
partition $X\in\parto$, if we have $f^{-1}(f(X))=X$, otherwise, we
say that it {\bf disregards} the partition $X$. If
$f^{-1}(f(X))=\{\omega\}$, then we say that $f$ {\bf completely
disregards} the partition $X$.

\begin{lm}\label{lm:respect}
For any surjection $f:\omega\surj\omega$ and for any $Z\in\parto$,
there is an $X\in\open{Z}$ such that $f$ either respects or
completely disregards the partition $X$.
\end{lm}

\proof For a surjection $f:\omega\surj\omega$, define the coloring
$\pi:\nseg{\omega}{2}\to\{0,1\}$ as follows. $\pi(s):=0$ if and only
if $f(s(0))\cap f(s(1))=\emptyset$. By Proposition~\ref{prop:thm},
there is a partition $X\in\open{Z}$ such that $\nseg{X}{2}$ is
monochromatic with respect to $\pi$, which implies that $f$ respects
$X$ in case of $\pi|_{\nseg{X}{2}}=\{0\}$, and $f$ completely
disregards $X$ is case of $\pi|_{\nseg{X}{2}}=\{1\}$. \eop

In the sequel we will use a slightly stronger version of
Proposition~\ref{prop:thm}, which is given in the following two
corollaries.

\begin{cor}\label{cor:weak-cor}
For any coloring of $\nseg{\omega}{(n+k+1)}$ with $r+1$ colors, where
$r,n,k\in\omega$, and for any dual-Ellentuck neighborhood
$\open{s,Y}$, where $|s|=n+1$, there is an infinite partition
$X\in\open{s,Y}$ such that $\nseg{s,X}{(n+k+1)}$ is monochromatic.
\end{cor}

\proof Let $\open{s,Y}$ be any dual-Ellentuck neighborhood, with
$|s|=n+1\ge 1$. Set $Y':=s\kap Y$, $R:=\bigcup_{i< n+1}Y'(i)$ and
$Y_R:=Y'\setminus\{Y'(i):i< n+1\}$, and take any order-preserving
bijection $f:\omega\setminus R\to \omega$. Then $Z:=f(Y_R)$ is an
infinite partition of $\omega$. For $u\in\nseg{Z}{n+k+1}$ we define
$\xi(u)\in\nseg{s,Y}{n+k+1}$ as follows.
$\mdom(\xi(u)):=f^{-1}(\mdom(u))$ and for $i<n+k+1$,
\begin{equation*} \xi(u)(i):=
                  \begin{cases}
                  \big( Y'(i)\cap\mdom(u)\big) \cup f^{-1}(u(i))
                  &\text{for $i<n+1$,}\\
                  f^{-1}(u(i)) &\text{otherwise.}
                  \end{cases}
\end{equation*}
Let $\pi:\nseg{\omega}{(n+k+1)}\to r+1$ be any coloring. Define
$\tau:\nseg{\omega}{(n+k+1)}\to r+1$ by stipulating $\tau (u):=\pi
(\xi(u))$. By Proposition~\ref{prop:thm} there is an infinite
partition $X'\in\open{Z}$ such that $\nseg{X'}{n+k+1}$ is
monochromatic with respect to the coloring $\tau$. Now let
$X\in\parto$ be such that
\begin{equation*} X(i):=
                 \begin{cases} Y'(i)\cup f^{-1}(X'(i))
                 &\text{for $i<n+1$}\\
                 f^{-1}(X'(i)) &\text{otherwise.}
                 \end{cases}
\end{equation*}
Then, by definition of $\tau$ and $X'$, $X\in\open{s,Y}$ and
$\nseg{s,X}{(n+k+1)}$ is monochromatic with respect to $\pi$. \eop

\begin{cor}\label{cor:strong-cor}
For any coloring of $\bigcup_{n\in\omega}\nseg{\omega}{(n+k+1)}$ with
$r+1$ colors, where $r,k\in\omega$, and for any $Z\in\parto$, there
is an infinite partition $X\in\open{Z}$ such that for any
$n\in\omega$ and for any $s\seg X$ with $|s|=n+1$,
$\nseg{s,X}{(n+k+1)}$ is monochromatic.
\end{cor}

\proof Using Corollary~\ref{cor:weak-cor} repeatedly, we can
construct the partition $X\in\parto$ straight forward by induction on
$n$. \eop

We say that a family $\nC\subs\parto$ has the {\bf
segment-coloring-property}, if for every coloring of
$\bigcup_{n\in\omega}\nseg{\omega}{(n+k+1)}$ with $r+1$ colors, where
$r,k\in\omega$, and for any $Z\in\nC$, there is an infinite partition
$X\in\open{Z}\cap\nC$, such that for any $n\in\omega$ and for any
$s\seg X$ with $|s|=n+1$, $\nseg{s,X}{(n+k+1)}$ is monochromatic.

If a partition-ultrafilter $\nU\in\puf{\parto}$ has the
segment-coloring-prop\-erty, then it is called a {\bf Ramseyan
ultrafilter}.

The next lemma shows that every partition-filter $\nF\in\pf{\parto}$
which has the segment-coloring-property is a partition-ultrafilter. A
similar result we have for Ramsey filters over $\omega$, since every
Ramsey filter is an ultrafilter.

\begin{lm}\label{lm:pf-puf}
If $\nF\subs\parto$ is a partition-filter which has the
segment-coloring-proper\-ty, then $\nF\subs\parto$ is a
partition-ultrafilter.
\end{lm}

\proof Take any $Z\in\parto$ such that for any $X\in\nF$, $Z\kap
X\in\parto$. Define the coloring $\pi:\nseg{\omega}{2}\to\{0,1\}$ by
stipulating $\pi(u)=0$ if and only if $u\in\nseg{Z}{2}$. Because
$\nF$ has the segment-coloring-property, there is a partition
$X\in\nF$ such that $\nseg{X}{2}$ is monochromatic with respect to
$\pi$, which implies that $X\ceq Z$ in case of
$\pi|_{\nseg{X}{2}}=\{0\}$, and $X\kap Z=\{\omega\}$ in case of
$\pi|_{\nseg{X}{2}}=\{1\}$. By the choice of $Z$ we must have $X\ceq
Z$, thus, since $\nF$ is a partition-filter, $Z\in\nF$. \eop

The following lemma gives a relation between Ramseyan and Ramsey
ultrafilters.

\begin{lm}\label{lm:RamseyRamseyan}
If $\nU$ is a Ramseyan ultrafilter, then
$\{\mMin(X)\setminus\{0\}:X\in\nU\}$ is a Ramsey ultrafilter over
$\omega$ (to be pedantic, one should say ``over
$\omega\setminus\{0\}$'').
\end{lm}

\proof Let $\tau:[\omega]^n\to r$ be any coloring of the $n$-element
subsets of $\omega$ with $r$ colors, where $n$ and $r$ are positive
natural numbers. Define $\pi:\nseg{\omega}{n}\to r$ by stipulating
$\pi(s):=\tau(\mMin(s^*)\setminus\{0\})$. Take $X\in\nU$ such that
$\nseg{X}{n}$ is monochromatic with respect to $\pi$, then, by the
definition of $\pi$, the set $[\mMin(X)\setminus\{0\}]^n$ is
monochromatic with respect to $\tau$. \eop

Ramsey ultrafilters over $\omega$ build the minimal points of the
Rudin-Keisler ordering on $\bNm$. This fact can also be expressed by
saying that a non-principal ultrafilter $\cU$ is a Ramsey ultrafilter
if and only if any function $g:\omega\to\omega$ is either constant or
one-to-one on some set of $\cU$. By Lemma~\ref{lm:respect}, we get a
similar result for Ramseyan ultrafilters with respect to the ordering
``$\les$''.

\begin{thm}\label{thm:Ruf-respect}
If $\nU$ is a Ramseyan ultrafilter, then for any surjection
$f:\omega\surj\omega$ there is an $X\in\nU$ such that $f$ either
respects or completely disregards $X$.
\end{thm}

\proof The proof is the same as the proof of Lemma~\ref{lm:respect},
but restricted to the partition-ultrafilter $\nU$. \eop

\subsection{On the existence of Ramseyan ultrafilters}\label{sec:ex}

As we have seen in Lemma~\ref{lm:RamseyRamseyan}, every Ramseyan
ultrafilter induces a Ramsey ultrafilter over $\omega$. It is not
clear if the converse holds as well. However, Ramseyan ultrafilters
are always forceable: Let $\UUs$ be the forcing notion consisting of
infinite partitions, stipulating $X\leq Y\Leftrightarrow X\ceq^* Y$.
$\UUs$ is the natural dualization of the forcing notion
$\la\cP(\omega)/{\rm fin},\subs^*\ra$, in the sequel denoted by
$\UU$, and it is not hard to see that if $\nG$ is $\UUs$-generic over
${\bfV}$, then $\nG$ is a Ramseyan ultrafilter in ${\bfV}[\nG]$.
Since $\UUs$ is $\sigma$-closed, as a consequence we get that
Ramseyan ultrafilters exist if we assume the continuum hypothesis
(denoted by $\CH$). On the other hand we know by
Lemma~\ref{lm:RamseyRamseyan} that Ramseyan ultrafilters cannot exist
if there are no Ramsey ultrafilters. Kenneth Kunen proved
({\cf}\cite[Theorem\,91]{Jechbook}) that it is consistent with $\ZFC$
that Ramsey ultrafilters don't exist. We like to mention that Saharon
Shelah showed that even $p$-points, which are weaker ultrafilters
than Ramsey ultrafilters, may not exist (see \cite[VI\,$\S
4$]{Shelah}). He also proved that it is possible that---up to
isomorphisms---there exists a unique Ramsey ultrafilter (see
\cite[VI\,$\S 5$]{Shelah}).

In the following, $\fc$ denotes the cardinality of the continuum and
$\twotoc$ denotes the cardinality of its power-set.

Andreas Blass proved that Martin's Axiom, denoted by $\MA$, implies
the existence of $\twotoc$ Ramsey ultrafilters (see
\cite[Theorem\,2]{BlassTrans}). He mentions in this paper that with
$\CH$ in place of $\MA$, this result is due to Keisler and with $1$
in place of $\twotoc$, it is due to Booth
({\cf}\cite[Theorem\,4.14]{Booth}). Further he mentions that his
proof is essentially the union of Keisler's and Booth's proof.
However, Blass' proof uses at a crucial point that $\MA$ implies that
the tower number is equal to $\fc$. Such a result we don't have for
partitions, because Timothy Carlson proved that the dual-tower number
is equal to $\aleph_1$ (see \cite[Proposition\,4.3]{Matet}). So,
concerning the existence of Ramseyan ultrafilters under $\MA$, we
cannot simply translate the proof of Blass, and it seems that $\MA$
and sets of partitions are quite unrelated. But as mentioned above,
if one assumes $\CH$, then Ramseyan ultrafilters exist. Moreover,
with respect to the equivalence relation ``$\simeq$'' (defined in
section~\ref{sec:ordering}) we get the following (for a similar
result w.r.t.\;the Rudin-Keisler ordering see
\cite[p.\,149]{BlassTrans}).

\begin{thm}\label{thm:CHmany}
$\CH$ implies the existence of $\twotoc$ pairwise non-equivalent
Ramseyan ultrafilters.
\end{thm}

\proof Assume ${\bfV}\models\CH$. Let $\chi$ be large enough such
that $\cP(\parto)\in H(\chi)$, {\ie}, the power set of $\parto$ (in
${\bfV}$) is hereditarily of size $<\chi$. Let $\bfN$ be an
elementary submodel of $\langle H(\chi),\in\rangle$ with $|\bfN
|={\aleph_1}$, containing all reals (or equivalently, all partitions)
of ${\bfV}$. We consider the forcing notion $\UUs$ in the model
$\bfN$. Since $|\bfN|=\aleph_1$, in ${\bfV}$ there is an enumeration
$\{D_\alpha\subs\parto: \alpha<{\omega_1}\}$ of all dense sets of
$\UUs$ which lie in $\bfN$. For any $Z\in\parto\cap\bfV$, let
$Y^{\alpha,0}_Z,Y^{\alpha,1}_Z\in D_\alpha$ be such that
$Y^{\alpha,0}_Z\ceq ^* Z$, $Y^{\alpha,1}_Z \ceq ^* Z$ and
$Y^{\alpha,0}_Z\kap Y^{\alpha,1}_Z\notin\parto$ (since $D_\alpha$ is
dense, such partitions exist). For any function $\zeta:\fc\to\{0,1\}$
we can construct a set $H_\zeta =\{X_\alpha: \alpha<{\omega_1}\}$ in
${\bfV}$ such that for all $\beta < \alpha<\omega_1$ we have
$X_\alpha\ceq^* Y^{\beta,\zeta(\beta)}_{X_\beta}$. By construction,
for any function $\zeta$, the set $G_\zeta:=\{X\in\parto:
X_\alpha\ceq^* X\ \text{for some $X_\alpha\in H_\zeta$}\}$ is
$\UUs$-generic over $\bfN$, thus, a Ramseyan ultrafilter in
$\bfN[G_\zeta]$, and since $\UUs$ is $\sigma$-closed and therefore
adds no new reals, $G_\zeta$ is also a Ramseyan ultrafilter in
$\bfV$. Furthermore, if $\zeta \neq \zeta'$, then the two Ramseyan
ultrafilters $G_\zeta$ and $G_{\zeta'}$ are different (consider the
two partitions $X_{\beta+1}\in H_\zeta$ and $X'_{\beta+1}\in
H_{\zeta'}$, where $\zeta(\beta)\neq\zeta'(\beta)$). Hence, in
$\bfV$, there are $\twotoc$ Ramseyan ultrafilters. Because there are
only $\fc$ surjections from $\omega$ onto $\omega$, no equivalence
class  (w.r.t.\;``$\simeq$'') can contain more than $\fc$ Ramseyan
ultrafilters, so, in $\bfV$, there must be $\twotoc$ pairwise
non-equivalent Ramseyan ultrafilters. \eop

\section{The happy families' relatives}

\subsection{Relatively happy families}

As we will see below, the partition-families which have the
segment-coloring-property are related to special coideals, so-called
happy families, which are introduced and rigorously investigated by
Adrian Mathias in \cite{Mathias}. So, partition-families with the
segment-coloring-property can be considered as ``relatives of happy
families''.

Let us first consider the definition of Mathias' happy families.

Let $\reals$ be the set of all infinite subsets of $\omega$, and let
$\fin$ be the set of all finite subsets of $\omega$. A set
$\cI\subs\cP(\omega)$ is a {\bf free ideal}, if $\cI$ is an ideal
which contains the {\bf Fr\'echet ideal} $\fin$. A set
$\cF\subs\cP(\omega)$ is a {\bf free filter}, if $\{y:\omega\setminus
y\in\cF\}$ is an ideal containing the Fr\'echet ideal. For
$a\in\fin$, let $a^*:=\mmax\{n+1:n\in a\}$, in particular, $0^*=0$.
For $x,y\in\cP(\omega)$ we write $y\subs^* x$ if $(y\setminus
x)\in\fin$. For a set $\cB\subs\cP(\omega)$, let $\fil(\cB)$ be the
free filter generated by $\cB$, so, $x\in\fil(\cB)$ if and only if
there is a finite set $y_0,\ldots,y_n\in\cB$ such that
$(y_0\cap\ldots\cap y_n)\subs^* x$.

A set $x\subs\omega$ is said to {\bf diagonalize} the family
$\{x_a:a\in\fin\}$, if $x\subs x_0$ and for all $a\in\fin$, if $\mmax
(a)\in x$, then $(x\setminus a^*)\subs x_a$.

The family $\cA\subs\cP(\omega)$ is {\bf happy}, if
$\cP(\omega)\setminus\cA$ is a free ideal and whenever
$\fil\{x_a:a\in\fin\}\subs\cA$, there is an $x\in\cA$ which
diagonalizes $\{x_a:a\in\fin\}$.

In terms of happy families one can define Ramsey ultrafilters as
follows: A Ramsey ultrafilter is an ultrafilter that is also a happy
family.

Now we turn back to partitions. The Fr\'echet ideal corresponds to
the set of finite partitions, and therefore, the notion of a free
filter corresponds to partition-filters containing only infinite
partitions, hence, to partition-filters $\nF\subs\parto$. For a set
$\nB\subs\parto$, let $\fil(\nB)$ be the partition-filter generated
by $\nB$, so, $X\in\fil(\nB)$ if and only if there is a finite set of
partitions $Y_0,\ldots,Y_n\in\nB$ such that $(Y_0\kap\ldots\kap
Y_n)\ceq^* X$.

A partition $X$ is said to diagonalize the family $\{X_s:s\in\NN\}$,
if $X\ceq X_{\emptyset}$ and for all $s\in\NN$, if $s^*\seg X$, then
$\big(\bigcup s^* \kap X\big) \ceq X_s$.

The family $\nA\subs\parto$ is {\bf relatively happy}, if whenever
$\fil\{X_s:s\in\NN\}\subs\nA$, there is an $X\in\nA$ which
diagonalizes $\{X_s:s\in\NN\}$.

An example of a relatively happy family is $\parto$, the set of all
infinite partitions (compare with \cite[Example\,0.2]{Mathias}).
Another example of a much smaller relatively happy family is given in
the following theorem (compare with \cite[p.\,63]{Mathias}).

\begin{thm}\label{thm:Ramseyanhappy}
Every Ramseyan ultrafilter is relatively happy.
\end{thm}

\proof Let $\nU\subs\parto$ be a partition-ultrafilter which has the
segment-coloring-property and let $\{X_s:s\in\NN\}\subs\nU$ be any
family. Since $\nU$ is a partition-filter, we obviously have
$\fil\{X_s:s\in\NN\}\subs\nU$. For $t\in\NN$ with $|t|\ge 2$, let
$s_t$ be such that $s_t^*\seg t$ and $|s_t|=|t|-2$. Define the
coloring $\pi:\bigcup_{n\in\omega}\nseg{\omega}{(n+2)}\to\{0,1\}$ by
stipulating
\begin{equation*} \pi(t):=\begin{cases} 0 &
                          \text{if $\big(\bigcup s_t^*\kap t^*\big)
                                        \ceq X_{s_t}$,}\\
                                        1 & \text{otherwise.}
                          \end{cases}
\end{equation*}
Let $X\in\open{X_{\emptyset}}\cap\nU$ be such that for any
$n\in\omega$ and for any $s^*\seg X$ with $|s|=n$,
$\nseg{s^*,X}{(n+2)}$ is monochromatic with respect to $\pi$. Take
any $s^*\seg X$. Since $\nseg{s^*,X}{(|s|+2)}$ is monochromatic with
respect to $\pi$, each $t^*\ceq X$ with $s^*\seg t$ and $|t|=|s|+2$
gets the same color. Hence, for all such $t$'s we have either
$\big(\bigcup s^*\kap t^*\big)\ceq X_s$, which implies $X\ceq^* X_s$,
or $\big(\bigcup s^*\kap t^*\big)\not\ceq X_s$, which implies $X\kap
X_s\notin\parto$. The latter is impossible, since it contradicts the
assumption that $\nU$ is a partition-filter. So, we are always in the
former case, which completes the proof. \eop

\subsection{A game characterization}

There is a characterization of happy ultrafilters over $\omega$,
{\ie}, of Ramsey ultrafilters, in terms of games
({\cf}\cite[Theorem\,4.5.3]{BartoszynskiJudah}). A similar
characterization we get for relatively happy partition-ultrafilter.

Let $\nU$ be a partition-ultrafilter. Define a game $G(\nU)$ played
by players~I and II as follows:

$$\begin{array}{ccccccccc}
 {\operatorname{I}} &\ \ \ &X_1 &    &X_2 &    &X_3 &    & \\
                    &      &    &    &    &    &    &    & \ldots \\
 {\operatorname{II}}&      &    &s_1 &    &s_2 &    &s_3 &
\end{array}$$

\noindent Player~I on the $n$-th move plays a partition $X_n\in\nU$.
Player~II responds with a segment $s_n\in\NN$ such that $|s_n|=n$,
$s_{n-1}^*\seg s_n$ and for all $m<n$, $\big(\bigcup s_m^*\kap
s_n^*\big)\ceq X_{m+1}$, where $s_0:=\emptyset$. Player~I wins if and
only if the unique partition $X$ with $s_n\seg X$ (for all $n$) is
not in $\nU$.

\begin{thm}\label{thm:game}
Let $\nU\in\puf{\parto}$, then player~I has a winning strategy in
$G(\nU)$ if and only if $\nU$ is not relatively happy.
\end{thm}

\proof Assume first that the partition-ultrafilter $\nU$ is
relatively happy and that $\{X_s:s\in\NN\}$ is a strategy for
player~I. This means, player~I begins with $X_{\emptyset}$ and then,
if $s_n$ is the $n$-th move of player~II, player~I plays $X_{s_n}$.
Because $\nU$ is relatively happy, there is a partition $X\in\nU$
which diagonalizes the family $\{X_s:s\in\NN\}$, in particular,
$X\ceq X_{\emptyset}$. Now, by the definition of $X$ and by the rules
of the game $G(\nU)$, player~II can play the segments of $X$. More
precisely, player~II plays on the $n$-th move the segment $s_n$, so
that $|s_n|=n$ and $s_n^*\seg X$. Since $X\in\nU$, the strategy
$\{X_s:s\in\NN\}$ was not a winning strategy for player~I.\\[1ex] Now
assume that the strategy $\sigma =\{X_s:s\in\NN\}$ is not a winning
strategy for player~I. Consider the game where player~I is playing
according to the strategy $\sigma$. In this game, player~II can play
segments $s_n$ such that the unique partition $X$ with $s_n\seg X$
(for all $n$) is in $\nU$. We have to show that $X$ diagonalizes the
family $\{X_s:s\in\NN\}$. For $n\in\omega$, let $s_n\in\NN$ be such
that $s_n^*\seg X$ and $|s_n|=n$. Fix $m\in\omega$, then, by the
rules of the game, for any $n>m$ we have $\big(\bigcup s_m^*\kap
s_n^*\big)\ceq X_{m+1}$, which implies $\big(\bigcup s_m^*\kap
X\big)\ceq X_{m+1}$. Since player~I follows the strategy $\sigma$,
$X_{m+1}=X_{s_m}$, and because $m$ was arbitrary, for all
$m\in\omega$ we get $\big(\bigcup s_m^*\kap X\big)\ceq X_{s_m}$.
Hence, $X$ diagonalizes the family $\{X_s:s\in\NN\}$. \eop

\section{The combinatorics of dual Mathias forcing}

First we recall the Ellentuck topology on $\reals$. For $x\in\reals$
and $a\in [\omega]^{<\omega}$ with $x\cap(\mmax(a)+1)=a$, let
$[a,x]^\omega:=\{y\in\reals : a\subs y\subs x\}$, and let the basic
open sets on $\reals$ be the sets $[a,x]^\omega$. These sets are
called {\bf Ellentuck neighborhoods}. The topology induced by the
Ellentuck neighborhoods is called {\bf Ellentuck topology}
({\cf}\cite{Ellentuck}).

The {\bf Mathias forcing} $\MM$, introduced in \cite{Mathias},
consists of ordered pairs $\la a,x\ra$ such that $[a,x]^\omega$ is an
Ellentuck neighborhood and the ordering on $\MM$ is defined by
stipulating $\la a,x\ra\leq\la b,y\ra
\;\Leftrightarrow\;[a,x]^\omega\subs [b,y]^\omega$.

Mathias forcing restricted to a non-principal ultrafilter $\cU$,
denoted by $\MMU{\cU}$, consists of the ordered pairs $\la
a,x\ra\in\MM$, where in addition we require that $x\in\cU$.

Mathias forcing has a lot of nice combinatorial properties (some of
them are mentioned below) which also hold for Mathias forcing
restricted to a Ramsey ultrafilter (see \cite{Mathias}).

The {\bf dual Ellentuck topology} on $\parto$ is the topology induced
by the dual Ellentuck neighborhoods (defined in
section~\ref{sec:coloring}). Now, the {\bf dual Mathias forcing}
$\MMs$, introduced in \cite{CarlsonSimpson}, is defined similarly to
Mathias forcing $\MM$, using the dual Ellentuck topology instead of
the Ellentuck topology. So, $\MMs$ consists of ordered pairs $\la
s,X\ra$ such that $\open{s,X}$ is a dual Ellentuck neighborhood and
the ordering on $\MMs$ is defined by stipulating $\la s,X\ra\leq\la
t,Y\ra\;\Leftrightarrow\;\open{s,X}\subs \open{t,Y}$.

Dual Mathias forcing restricted to a partition-ultrafilter
$\nU\in\puf{\parto}$, denoted by $\MMUs{\!\nU}$, consists of the
ordered pairs $\la s,X\ra\in\MMs$, where in addition we require that
$X\in\nU$ (see {\eg}\;\cite{Lorisym} and \cite{boise}).

Both, Mathias forcing as well as dual Mathias forcing, are proper
forcings. Moreover, both have (i)\;a decomposition, (ii)\;pure
decision and (iii)\;the homogeneity property (see
{\eg}\;\cite{Mathias}, \cite{CarlsonSimpson} and \cite{Lorisym}):
\begin{enumerate}
\item[{(i)}] {\bf Decomposition}: $\MM\approx\UU *\MMU{\bfU}$, where
$\bfU$ is the canonical $\UU$-name for the $\UU$-generic object
($\UU$ as in section~\ref{sec:ex}).\smallskip

$\MMs\approx\UUs *\MMUs{\bfU}$, where $\bfU$ is the canonical
$\UUs$-name for the $\UUs$-generic object ($\UUs$ as in
section~\ref{sec:ex}).

\item[{(ii)}] {\bf Pure decision}: For any $\MM$-condition $\la a,x\ra$
and any sentence $\Phi$ of the forcing language $\MM$, there is an
$\MM$-condition $\la a,y\ra \leq \la a,x\ra$ such that either $\la
a,y\ra \forces{\MM}\Phi$ or $\la
a,y\ra\forces{\MM}\neg\Phi$.\smallskip

For any $\MMs$-condition $\la s,X\ra$ and any sentence $\Phi$ of the
forcing language $\MMs$, there is an $\MMs$-condition $\la s,Y\ra
\leq \la s,X\ra$ such that either $\la s,Y\ra \forces{\MMs}\Phi$ or
$\la s,Y\ra\forces{\MMs}\neg\Phi$.

\item[{(iii)}] {\bf Homogeneity property}: If $x_G$ is $\MM$-generic
over $\bfV$ and $y\in [x_G]^\omega$, then $y$ is also $\MM$-generic
over $\bfV$.\smallskip

If $X_G$ is $\MMs$-generic over $\bfV$ and $Y\in\open{X_G}$, then $Y$
is also $\MMs$-generic over $\bfV$.
\end{enumerate}

In \cite{Lorisym} it is shown that if $\nF\subs\parto$ is a so-called
{\em game-family}, then $\MMUs{\nF}$ has pure decision and the
homogeneity property (\cite[Thm.\,4.3\,\&\,4.4]{Lorisym}).
Game-families have the segment-coloring-property and therefore, the
so-called {\em game-filters}, {\ie}, game-families which are
partition-filters, are Ramseyan ultrafilters. Unlike for Ramseyan
ultrafilters, it is not clear if $\CH$ implies the existence of
game-filters, so, it seems that game-filters are stronger than
Ramseyan ultrafilters. However, in the sequel we show that if
$\nU\in\puf{\parto}$ is a Ramseyan ultrafilter, then $\MMUs{\!\nU}$
has pure decision and the homogeneity property.

Recently, Stevo Todor\v{c}evi\'c gave an abstract presentation of
Ellentuck's theorem by introducing the notion of a {\it quasi
ordering with approximations\/} which {\it admits a finitization\/}
and the notion of a {\it Ramsey space}. The {\sc Abstract Ellentuck
Theorem} says that a quasi ordering with approximations which admits
a finitization and satisfies certain axioms is a Ramsey space.

Let $\nU\in\puf{\parto}$ be a partition-ultrafilter and let
``$\ceq$'' be the quasi ordering on $\nU$. For each $n\in\omega$, let
the function $p_n: \nU\to\NN$ be such that $p_n(X)$ is the unique $s$
with $s^*\seg X$ and $|s|=n$. Let $p$ be the sequence
$(p_n)_{n\in\omega}$. It is easy to verify that the triple $(\nU,
\ceq,p)$ is a {\bf quasi ordering with approximations}. For
$n,m\in\omega$ and $X,Y\in\nU$ define:
$p_n(X)\ceq_{\text{fin}}p_m(Y)$ if and only if
$\mdom\big(p_n(X)\big)=\mdom\big(p_m(Y)\big)$ and $p_n(X)\ceq
p_m(Y)$. This definition verifies that $(\nU,\ceq,p)$ {\bf admits a
finitization}. If $\open{s,X}$ is a dual Ellentuck neighborhood and
$X\in\nU$, then $\open{s,X}\cap\nU$ is called a {\bf $\nU$-dual
Ellentuck neighborhood}. The topology on $\nU$, induced by the
$\nU$-dual Ellentuck neighborhoods, is called the {\bf $\nU$-dual
Ellentuck topology}. With respect to the $\nU$-dual Ellentuck
topology, the topological space $\nU$ is a {\bf Ramsey space}, if for
any subset $S\subs\nU$ which has the Baire property with respect to
the $\nU$-dual Ellentuck topology, and for any $\nU$-dual Ellentuck
neighborhood $\open{s,Y}\cap\nU$, there is a partition
$X\in\open{s,Y}\cap\nU$ such that either $\open{s,X}\cap\nU\subs S$
or $\open{s,X}\cap\nU\subs \nU\setminus S$.

Let $\nU\in\puf{\parto}$ be a Ramseyan ultrafilter. Since the triple
$(\nU,\ceq,p)$ satisfies certain axioms, by Todor\v{c}evi\'c's {\sc
Abstract Ellentuck Theorem}, the Ramseyan ultrafilter $\nU$ with
respect to the $\nU$-dual Ellentuck topology is a Ramsey space.
Moreover, we get the following two results.

\begin{thmm}\label{thm:pure}
If $\nU$ is a Ramseyan ultrafilter, then $\MMUs{\!\nU}$ has pure
decision.
\end{thmm}

\proof Let $\Phi$ be any sentence of the forcing language
$\MMUs{\!\nU}$. With respect to $\Phi$ we define $$D_0:=\{Y\in\nU:
\text{ for some $t\seg Y$, }\la t,Y\ra\forces{\MMUs{\!\nU}}
\neg\Phi\}\,,$$ and $$D_1:=\{Y\in\nU:\text{ for some $t\seg Y$, }\la
t,Y\ra\forces{\MMUs{\!\nU}} \Phi\}\,.$$ Clearly $D_0$ and $D_1$ are
both open (w.r.t.\;the $\nU$-dual Ellentuck topology) and $D_0\cup
D_1$ is dense (w.r.t.\;the partial order in $\MMUs{\!\nU}$). Because
$\nU$ is a Ramsey space, for any $\nU$-dual Ellentuck neighborhood
$\open{s,Y}\cap\nU$ there is an $X\in\open{s,Y}\cap\nU$ such that
$\open{s,X}\cap\nU\subs D_0$ or $\open{s,X}\cap\nU\cap D_0
=\emptyset$. In the former case we have $\la s,X\ra
\forces{\MMUs{\!\nU}}\neg\Phi$ and we are done. In the latter case we
find $X'\in\open{s,X}\cap\nU$ such that $\open{s,X'}\cap\nU\subs
D_1$. (Otherwise we would have $\open{s,X'}\cap\nU\cap (D_0\cup D_1)
=\emptyset$, which is impossible by the density of $D_0\cup D_1$.)
Hence, $\la s,X'\ra\forces{\MMUs{\!\nU}}\Phi$. \eop

\begin{thmm}\label{thm:homogen}
If $\nU$ is a Ramseyan ultrafilter, then $\MMUs{\!\nU}$ has the
homogeneity property.
\end{thmm}

\proof For a dense set $D\subs\MMUs{\!\nU}$, let $$\bigcup
D:=\{X\in\parto: X\in\open{s,Y}\text{ for some $\la s,Y\ra\in
D$}\}\,.$$ It is clear that a partition $X_G$ is
$\MMUs{\!\nU}$-generic if and only if $X_G\in\bigcup D$ for each
dense set $D\subs\MMUs{\!\nU}$. Let $D\subs\MMUs{\!\nU}$ be an
arbitrary dense set and let $D'$ be the set of all $\la
s,Z\ra\in\MMUs{\!\nU}$ such that $\open{t,Z}\subs \bigcup D$ for all
$t\ceq s$ with $\mdom (t)=\mdom (s)$.\\[1ex] First we show that $D'$
is dense in $\MMUs{\!\nU}$. For this, take an arbitrary $\la
s,W\ra\in\MMUs{\!\nU}$ and let $\{t_i: 0\leq i\leq m\}$ be an
enumeration of all $t\in\NN$ such that $t\ceq s$ and $\mdom (t)=\mdom
(s)$. Because $D$ is dense in $\MMUs{\!\nU}$, $\bigcup D$ is open
(w.r.t.\;the $\nU$-dual Ellentuck topology), and since $\nU$ is a
Ramsey space, for every $t_i$ we find a $W'\in\nU$ such that $t_i\ceq
W'$ and $\open{t_i, W'}\subs\bigcup D$. Moreover, if we define
$W_{-1}:=W$, for every $i\leq m$ we can choose a partition
$W_i\in\nU$ such that $W_i\ceq W_{i-1}$, $s\seg W_i$ and
$\open{t_i,W_i}\subs\bigcup D$. Thus, $\la s,W_m\ra\in D'$, and
because $\la s,W_m\ra\le\la s,W\ra$, $D'$ is dense in
$\MMUs{\!\nU}$.\\[1ex] Let $X_G$ be $\MMUs{\!\nU}$-generic and let
$Y\in\open{X_G}$ be arbitrary. Since $D'$ is dense, there is a
condition $\la s,Z\ra\in D'$ such that $s\seg X_G\ceq Z$. Since
$Y\in\open{X_G}$, we have $t\seg Y\ceq Z$ for some $t\ceq s$ with
$\mdom(t)=\mdom(s)$, and because $\open{t,Z}\subs\bigcup D$, we get
$Y\in\bigcup D$. Hence, $Y\in\bigcup D$ for each dense set
$D\subs\MMUs{\!\nU}$, which completes the proof. \eop

\section*{Appendix}

In this section we are gathering some results concerning the dual
form of some cardinal characteristics of the continuum. For the
definition of the classical cardinal characteristics, as well as for
the relation between them, we refer the reader to \cite{Vaughan}.

First we consider the shattering cardinal $\fh$. This cardinal was
introduced in \cite{Balcar.etal} as the minimal height of a tree
$\pi$-base of $\bNm$. Later it was shown by Szymon Plewik in
(\cite{Plewik}) that $\fh=\add(r^0)=\cov(r^0)$, where $r^0$ denotes
the ideal of Ramsey-null sets. It is easy to see that $\fp\le\fh$,
and therefore, $\MA(\sigma\text{-centered})$ implies $\fh=\fc$.

The dual form of the classical cardinal characteristics were
introduced and investigated in \cite{Cichon.etal} and further
investigated in \cite{Lorisha}. Concerning the dual-shattering
cardinal $\fH$, one easily gets $\aleph_1\le\fH\le\fh$, and in
\cite{Lorisha} it is shown that $\fH >\aleph_1$ is consistent
relative to $\ZFC$ and that $\fH=\add(R^0)=\cov(R^0)$, where $R^0$
denotes the ideal of dual Ramsey-null sets. After all these
symmetries, one would not expect the following: $\MA+ (\fc >\fH)$ is
consistent relative to $\ZFC$. This was proved by J\"org Brendle in
\cite{Brendle} and implies that $\fH <\fp$ is consistent relative to
$\ZFC$.

Concerning the reaping and the dual-reaping number $\fr$ and $\fR$,
respectively, the situation looks different. It is shown in
\cite{Lorisym} that $\fp\le\fR\le\mmin\{\fr,\fri\}$, and thus we get
$\MA(\sigma\text{-centered})$ implies $\fR=\fc$. Further, it is easy
to show that $\fR\le\fU$, where $\fU$ denotes the
partition-ultrafilter base number, {\ie}, the dual form of $\fu$, and
consequently, $\MA(\sigma\text{-centered})$ implies $\fU=\fc$.

For a Ramsey ultrafilter $\cU$, Brendle introduced in
\cite{BrendleFund} the ideal $r_{\cU}^0$, which is the ideal of
Ramsey-null sets with respect to the ultrafilter $\cU$. Concerning
this ideal $r_{\cU}^0$, he showed for example that
$\fhom\le\non(r_{\cU}^0)$, where $\fhom$ is the homogeneity number
investigated by Blass in \cite[Section\,6]{BlassSimple}. There, Blass
also investigated the so-called partition number $\fpar$ and showed
that $\fpar = \mmin\{\fb,\fs\}$. Now, replacing the Ramsey
ultrafilter $\cU$ by a Ramseyan ultrafilter $\nU$, one obtains the
ideal $R^0_{\!\nU}$ of dual Ramsey-null sets with respect to $\nU$ as
the dualization of the ideal $r_{\cU}^0$, and replacing the colorings
of $[\omega]^2$---involved in the definition of $\fhom$ and
$\fpar$---by colorings of $\nseg{\omega}{2}$, one obtains the
cardinal characteristics $\fHom$ and $\fPar$ and could begin to
investigate them. But this is left to the reader. \ding{166}


\begin{thebibliography}{M}
\bibitem[HL\"o$\infty_1$]{luxor}
{\toappear{Lorenz~Halbeisen and Benedikt~L\"owe}{Ultrafilter
           spaces on the semilattice of partitions}
          {Topology and its Applications}}
\end{thebibliography}
\end{document}